\documentclass[reqno, 12pt]{amsart}
\usepackage{a4wide}
\usepackage{amscd}

\newtheorem{thm}{Theorem}[section]

\newtheorem{prop}[thm]{Proposition}
\newtheorem{cor}[thm]{Corollary}

\font\fr=eufm10 scaled \magstep1
\font\ffr=eufm8 scaled \magstep1

\def\C{\mathbb C}

\def\F{\mathbb F}

\def\H{\mathbb H}

\def\O{\mathbb O}

\def\R{\mathbb R}

\def\Z{\mathbb Z}

\def\ca{\hbox{\fr a}}

\def\cg{\hbox{\fr g}}
\def\ch{\hbox{\fr h}}

\def\ck{\hbox{\fr k}}

\def\cn{\hbox{\fr n}}

\def\cp{\hbox{\fr p}}

\def\cs{\hbox{\fr s}}
\def\ct{\hbox{\fr t}}

\def\cv{\hbox{\fr v}}
\def\cw{\hbox{\fr w}}

\def\cz{\hbox{\fr z}}

\def\cJ{\hbox{\fr J}}

\def\cck{\hbox{\ffr k}}

\def\pt{\hbox{\rm pt}}
\def\ad{\hbox{\rm ad}}
\def\Ad{\hbox{\rm Ad}}

\def\Im{\hbox{\rm Im}}

\def\qed{\hfill \hbox{\rm q.e.d.}}
\def\Spin{\hbox{\rm Spin}}

\begin{document}

\begin{center}
{\large\bf Cohomogeneity one actions on noncompact\\ symmetric spaces of rank one}\\
\bigskip\bigskip
By {\it J\"{u}rgen Berndt} and {\it Hiroshi Tamaru}\\
\end{center}
\bigskip\bigskip
{\small {\bf Abstract.} We classify, up to orbit equivalence, all
cohomogeneity one actions on the hyperbolic planes over the
complex, quaternionic and Cayley numbers, and on the complex
hyperbolic spaces $\C H^n$, $n \geq 3$. For the quaternionic
hyperbolic spaces $\H H^n$, $n \geq 3$, we reduce the
classification problem to a problem in quaternionic linear algebra
and obtain partial results. For real hyperbolic spaces, this
classification problem was essentially solved by
\'{E}lie Cartan.}\\


\thispagestyle{empty}


\section{Introduction}

An isometric action on a Riemannian manifold is of cohomogeneity
one if its orbit space is one-dimensional. Cohomogeneity one
actions are of current interest for the construction of
geometrical structures on manifolds, e.g.\ Einstein metrics and
metrics with special holonomies. The reason is that a
cohomogeneity one action can be used to reduce the system of
partial differential equations describing such a geometrical
structure to a nonlinear ordinary differential equation for which
one might be able to find explicit solutions. Given a Riemannian
manifold $M$, it is natural to find all cohomogeneity one actions
on it, perhaps just up to orbit equivalence. Two cohomogeneity one
actions on $M$ are orbit equivalent if there exists an isometry of
$M$ that maps the orbits of one action onto the orbits of the
other action. It is worthwhile to mention that the classification
problem of cohomogeneity one actions up to orbit equivalence is
equivalent to the classification problem of homogeneous
hypersurfaces up to isometric congruence. The latter one is a
classical problem is submanifold theory.

The cohomogeneity one actions on the spheres, equipped with their standard metric
of constant curvature, have been classified by Hsiang and Lawson \cite{HL}.
Remarkably, any such action is orbit equivalent to the isotropy representation of a
Riemannian symmetric space of rank two. For the other compact
symmetric spaces of rank one the classifications were obtained by
Takagi \cite{Ta} for the complex projective spaces and by Iwata \cite{Iw1},
\cite{Iw2} for the quaternionic projective spaces and the Cayley projective plane.
For simply connected
irreducible Riemannian symmetric spaces of higher rank the classification
was established by Kollross \cite{Ko}.

The methods employed by the above authors do not work for the noncompact dual
symmetric spaces. The noncompactness of the isometry group turns out to be
a subtle point. There can be uncountably many families of nonisomorphic
subgroups of the isometry group that act orbit equivalently by cohomogeneity
one. By using the classification of isoparametric hypersurfaces
on the Euclidean space $\R^n$ and the real hyperbolic
space $\R H^n$ by
Levi Civita \cite{Le}, Segre \cite{Se} and Cartan \cite{Ca}, one can obtain all
cohomogeneity one actions on these spaces up to orbit equivalence. In both cases
the orbit structure is either a Riemannian foliation, or a totally geodesic
subspace together with the distance tubes around it.
It is a general fact that a cohomogeneity one action on a symmetric
space of noncompact type, or more general on a Hadamard manifold, induces
either a Riemannian foliation, or has exactly one singular orbit and the generic orbits
are the distance tubes around it (see \cite{BB} for details and references).

In \cite{BT1} we obtained the classification, up to orbit equivalence,
of all cohomogeneity one actions on irreducible symmetric spaces of noncompact type
that induce a Riemannian foliation, that is, have no singular orbit. A
surprising consequence of this result is that the moduli space of all such
actions just depends on the rank of the symmetric space and possible duality or
triality principles on the space. In particular, on each noncompact symmetric space
of rank one this moduli space consists just of two elements. The corresponding foliations
are the horosphere foliation and a foliation with exactly one minimal leaf whose
geometry has been investigated in \cite{Be}.

The classification of all cohomogeneity one actions on irreducible
symmetric spaces of noncompact type that have a totally geodesic
singular orbit has been achieved in \cite{BT2}. It thus remains to
investigate the case of a non-totally geodesic singular orbit. As
mentioned above, in case of $\R^n$ and $\R H^n$ a singular orbit
is necessarily totally geodesic. It is remarkable that this is no
longer true for the other noncompact symmetric spaces of rank one:
the complex hyperbolic spaces $\C H^n$ ($n \geq 2$), the
quaternionic hyperbolic spaces $\H H^n$ ($n \geq 2$), and the
Cayley hyperbolic plane $\O H^2$. The first author and Br\"{u}ck
constructed in \cite{BB} many examples of cohomogeneity one
actions on these hyperbolic spaces (except for $\C H^2$) with a
non-totally geodesic singular orbit. The main result of this paper
says that, up to orbit equivalence, there are no further
cohomogeneity one actions on $\C H^n$ ($n \geq 3$), $\H H^2$ and
$\O H^2$. We also show that every singular orbit of a
cohomogeneity one action on $\C H^2$ is totally geodesic. For the
quaternionic hyperbolic space $\H H^n$, $n \geq 3$, we prove that
the set of orbit equivalence classes of cohomogeneity one actions
with a singular orbit of codimension $2$ is parametrized by the
closed interval $[0,\pi/2]$.

The results of this paper were partially obtained during a
common visit to the Mathe\-matical Research Institute Oberwolfach (Research
in Pairs programme). We would like to thank the Institute for its support and
the hospitality.
The second author was partially supported by
Grant-in-Aid for Young Scientists (B) 14740049,
The Ministry of Education, Culture, Sports, Science and Technology,
Japan.

\section{Preliminaries}

Let $M$ be a noncompact symmetric space of rank one. Then $M$ is either a
real hyperbolic space $\R H^n$, a complex hyperbolic space $\C H^n$, a
quaternionic hyperbolic space $\H H^n$, or a Cayley hyperbolic plane $\O H^2$,
where $n \geq 2$. We denote by $\F$ one of the real division algebras
$\R$, $\C$, $\H$ or $\O$, and by $\F H^n$ the corresponding hyperbolic space,
where we assume $n = 2$ if $\F = \O$.
Let $G$ be the identity component of the full isometry group of $M$, that is,
$G = SO^o(n,1),SU(n,1),Sp(n,1),F_4^{-20}$ for $\F=\R,\C,\H,\O$, respectively.
We fix a point $o \in M$ and denote by $K$ the isotropy subgroup of $G$ at $o$,
that is, $K = SO(n),S(U(n)U(1)),Sp(n)Sp(1),\Spin(9)$. Then, as a homogeneous space,
$M$ is isomorphic to $G/K$. We denote by $\cg$ and $\ck$ the Lie algebra of $G$
and $K$. Let $B$ be the Killing form of $\cg$ and $\theta$ the Cartan involution of
$\cg$ with respect to $\ck$. Then $\langle X,Y \rangle = -B(X,\theta Y)$ is a positive
definite inner product on $\cg$. Let $\cg = \ck + \cp$ be the Cartan decomposition
of $\cg$ induced by $\theta$. The restriction of $\langle \cdot, \cdot \rangle$ to
$\cp$ induces a Riemannian metric on $G/K$ turning it into a Riemannian symmetric space
of rank one.
We normalize the Riemannian metric on $M$ so that it becomes isometric to $G/K$
with the induced metric we just described.

Let $\ca$ be a maximal abelian subspace of $\cp$, which is just a one-dimensional
linear subspace since the rank of $M$ is one, and
$$
\cg = \cg_{-2\alpha} + \cg_{\alpha} + \cg_0 + \cg_{\alpha} + \cg_{2\alpha}
$$
the corresponding restricted root space decomposition of $\cg$. Note that
$\cg_{-2\alpha}$ and $\cg_{2\alpha}$ are trivial if $\F = \R$. Then
$$
\cg = \ck + \ca + \cn\ \ \ {\rm with}\ \ \ \cn = \cg_{\alpha} + \cg_{2\alpha}
$$
is an Iwasawa decomposition of $\cg$. The subalgebra $\cn$ of $\cg$ is abelian
if $\F = \R$ and two-step nilpotent otherwise. In fact, $\cn$ is isomorphic to
the $(2n-1)$-dimensional Heisenberg algebra if $\F = \C$, and to a certain
generalized Heisenberg algebra if $\F \in \{\H,\O\}$
(see \cite{BTV} for more details on this). Moreover,
$\cz = \cg_{2\alpha}$ is the centre of $\cn$ and equal to the derived subalgebra
$[\cn,\cn]$ of $\cn$. The dimension of $\cz$ is equal to $1,3,7$ for $\F = \C,\H,\O$,
respectively. The subalgebra $\ca + \cn$ of $\cg$ is solvable and $\cn$ is the derived
subalgebra of $\ca + \cn$.

We denote by $A$ resp.\ $N$ the connected closed subgroup of $G$ with
Lie algebra $\ca$ resp.\ $\cn$. Then $G = KAN$ is an Iwasawa decomposition of $G$ and,
since $K$ is the isotropy subgroup of $G$ at $o$, the solvable subgroup $AN$ of $G$
acts simply transitively on $M$.
Thus $M$ is isometric to the solvable Lie group $AN$
equipped with a suitable left-invariant Riemannian metric.

We define $\cv = \cg_{\alpha}$. Then we can identify $\cv$ with $\R^{n-1},\C^{n-1},\H^{n-1},\O$
for $\F = \R,\C,\H,\O$, respectively. More precisely, if $\F = \R$ then $\cv$ is isomorphic
to $\R^{n-1}$ as a real vector space. If $\F = \C$, the K\"{a}hler structure on $\C H^n$ induces
a complex vector space structure on $\cv$ so that it becomes isomorphic to $\C^{n-1}$,
and if $\F = \H$, the quaternionic K\"{a}hler structure on $\H H^n$ induces a (right) quaternionic
vector space structure on $\cv$ so that it becomes isomorphic to the (right) quaternionic
vector space $\H^{n-1}$. For $\F = \O$ we simply identify $\O$ with $\R^8$, and $\cv$ becomes
isomorphic to $\O$ as a real vector space.

\section{The reduction}\label{reduction}

In this section we reduce our classification problem to the problem of classifying
certain subalgebras of $\ca + \cn$.
We start with a general result about cohomogeneity one actions on Hadamard manifolds, i.e.,
connected, simply connected, complete Riemannian manifolds
of nonpositive curvature. Every symmetric space of noncompact type is a Hadamard manifold.
We recall that Cartan's Fixed Point Theorem states that the action of any compact
subgroup of the isometry group of a Hadamard manifold has a fixed point. We refer to
\cite{Eb} for more details on Hadamard manifolds and Cartan's Fixed Point Theorem.

\begin{prop}\label{solvable}
Let $M$ be a Hadamard manifold and $H$ a connected subgroup of
the isometry group of $M$ that acts
with cohomogeneity one on $M$ and has a singular orbit $F$. Then
there exists a connected solvable subgroup of $H$ that acts transitively on $F$.
\end{prop}

{\it Proof.} We choose a Levi-Malcev decomposition
$$
\ch = \ch_{ss} + \ch_{solv}
$$
of the Lie algebra $\ch$ of $H$
into the semidirect sum of a semisimple subalgebra $\ch_{ss}$ and a solvable
ideal $\ch_{solv}$. For the semisimple subalgebra $\ch_{ss}$ we choose an Iwasawa decomposition
$$
\ch_{ss} = \ch_{ss}^{cpct} + \ch_{ss}^{solv}
$$
of $\ch_{ss}$ into the vector space direct sum of a compact subalgebra $\ch_{ss}^{cpct}$
and a solvable subalgebra $\ch_{ss}^{solv}$. Then
$$
\ch = \ch_{ss}^{cpct} + (\ch_{ss}^{solv} + \ch_{solv})\ ,
$$
where $\ch_{ss}^{solv} + \ch_{solv}$ is a semidirect sum of the two solvable subalgebras
so that $\ch_{solv}$ is the ideal in it. Note that the semidirect sum of two
solvable Lie algebras is solvable as well. We denote by $H_{ss}^{cpct}$ and
$H^{solv}$ the connected subgroup of $H$ with Lie algebra
$\ch_{ss}^{cpct}$ and $\ch_{ss}^{solv} + \ch_{solv}$, respectively.

By Cartan's Fixed Point Theorem, there exists a point $p \in M$ that is fixed
under the action of the compact group $H_{ss}^{cpct}$. If $p \in F$, then clearly the
solvable group $H^{solv}$ acts transitively on $F$. If $p \notin F$,
then $p$ is on a principal orbit of the $H$-action on $M$, and it follows that the
solvable group $H^{solv}$ acts transitively on this principal orbit.
Since the action of $H$ on $M$ is of cohomogeneity one, we easily see that
$H^{solv}$ acts transitively on each orbit of the $H$-action and,
in particular, also transitive on the singular orbit $F$.
This finishes the proof of Proposition \ref{solvable} \qed

\medskip
We denote by $M(\infty)$ the ideal boundary of $M$ whose points are given by
the equivalence classes of asymptotic geodesics in $M$, and equip
$\bar{M} = M \cup M(\infty)$ with the cone topology. The action of $H$ on $M$
extends canonically to an action of $H$ on $\bar{M}$.

From now on we assume that $M = \F H^n$ and that the singular orbit $F$
of the cohomogeneity one action by $H$ on $M$ is
not totally geodesic. Then none of the $H$-orbits on $M$ is totally geodesic,
and a result by Alekseevsky and Di Scala \cite{AD} implies that there exists
a unique point $x \in M(\infty)$ that is fixed under the $H$-action on $M(\infty)$.
We fix a point $o \in F$ and consider the Iwasawa decomposition
$$
\cg = \ck + \ca + \cn
$$
that is determined by $o$ and $x$. Since $H \cdot x = x$, we have
$$
\ch \subset \ck_x + \ca + \cn\ ,
$$
where $\ck_x$ is the centralizer of $\ca$ in $\ck$. We denote by $K_x$ the
connected subgroup of $K$ with Lie algebra $\ck_x$. Then we have
$K_x = SO(n-1),S(U(n-1)U(1)),Sp(n-1)Sp(1),\Spin(7)$ for $\F = \R,\C,\H,\O$,
respectively, and $H \subset K_xAN$.

By Proposition \ref{solvable} there exists a solvable subgroup $S$ of $H$ that
acts transitively on the singular orbit $F$. We denote by $\cs$ the subalgebra of
$\ch$ corresponding to $S \subset H$. We recall that we may write the nilpotent
subalgebra $\cn$ in the form $\cn = \cg_{\alpha} + \cg_{2\alpha}$ with some suitable
root spaces $\cg_{\alpha}$ and $\cg_{2\alpha}$. Since $\ck_x$ centralizes $\ca$,
it normalizes each root space and hence $\cn$, which implies that $\ca + \cn$
is an ideal in $\ck_x + \ca + \cn$. Thus the canonical projection
$$
\pi : \ck_x + \ca + \cn \to \ck_x
$$
is a Lie algebra homomorphism, and it follows that
$$
\cs_c = \pi(\cs)
$$
is a solvable subalgebra of $\ck_x$. Since every solvable subalgebra of
a compact Lie algebra is abelian, we conclude that
\begin{equation}\label{sc-abelian}
\cs_c\ {\rm is\ an\ abelian\ subalgebra\ of}\ \ck_x.
\end{equation}

Let
$$
\tau : \ck_x + \ca + \cn \to \ca + \cn
$$
be the canonical projection and define
$$
\cs_n = \tau(\cs)\ .
$$
It is clear that
\begin{equation}\label{dimensions}
\dim \cs_n = \dim F\ .
\end{equation}
Our aim is to show that $\cs_n$ is a subalgebra of $\ca + \cn$ and that the
orbit through $o$ of the action of the corresponding subgroup $S_n$ of $AN$ is just the
singular orbit $F$.
For each $k \in K_x$ the differential $d_ok$ of $k$ at $o$ is given by
$d_ok = \Ad(k)|(\ca+\cn)$, where we identify $T_oM$ with $\ca + \cn$
by means of $M = G/K = AN$. Since the isotropy subgroup $H_o$ of $H$ at $o$
acts transitively on the unit sphere in the normal space $\nu_oF$ of $F$ at $o$,
and as $\ck_x$ centralizes $\ca$, we necessarily have
\begin{equation}\label{normal-in-n}
\nu_oF \subset \cn
\end{equation}
and hence
\begin{equation}\label{a-in-sn}
\ca \subset \cs_n\ .
\end{equation}

We shall now prove that
\begin{equation}\label{derived-in-sn}
[\cs,\cs] = \cs_n \cap \cn \ .
\end{equation}
Since $\cs \subset \cs_c + \cs_n$, we have
$$
[\cs,\cs] \subset [\cs_c,\cs_c] + [\cs_c,\cs_n] + [\cs_n,\cs_n]\ .
$$
The subalgebra $[\cs_c,\cs_c]$ is trivial since $\cs_c$ is abelian according to (\ref{sc-abelian}).
Since $\cs_c \subset \ck_x$, $\cs_n \subset \ca + \cn$ and $\ck_x$
centralizes $\ca$ and normalizes $\cn$, we have
$[\cs_c,\cs_n] \subset [\ck_x,\ca + \cn] \subset \cn$.
Finally, since $\cs_n \subset \ca + \cn$ and $\cn$ is the derived subalgebra
of $\ca + \cn$, we see that $[\cs_n,\cs_n] \subset \cn$. Altogether this implies
$[\cs,\cs] \subset \cn$, which readily yields $[\cs,\cs] \subset \cs_n \cap \cn$.
For the converse, we fix the element $B \in \ca$ for which $[B,V] = V$ and $[B,Z] = 2Z$
holds for all $V \in \cv = \cg_{\alpha}$ and $Z \in \cz = \cg_{2\alpha}$.
Because of (\ref{a-in-sn}) there exists an elements $\tilde{B} \in \cs_c$
so that $\tilde{B} + B \in \cs$. Let $X = V + Z \in \cv + \cz$ be an arbitrary
element in the orthogonal complement of $[\cs,\cs]$ in $\cs_n \cap \cn$. Then there exists
a vector $\tilde{X} \in \cs_c$ so that $\tilde{X} + X \in \cs$, and we have
$$
0 = \langle X , [\tilde{B} + B, \tilde{X} + X] \rangle
= \langle X, [\tilde{B},\tilde{X}] + [\tilde{B},X] + [B,\tilde{X}] + [B,X] \rangle\ .
$$
Since $\cs_c$ is abelian we have $[\tilde{B},\tilde{X}] = 0$. Since ${\rm ad}(\tilde{B})$ is
a skewsymmetric transformation we have $\langle X,[\tilde{B},X] \rangle = 0$. And since
$\cs_c \subset \ck_x$ and $\ck_x$ centralizes $\ca$ we have $[B,\tilde{X}] = 0$. This implies
$$
0 = \langle X,[B,X] \rangle = \langle V+Z , V+2Z \rangle =
\langle V,V \rangle + 2 \langle Z,Z \rangle\ ,
$$
and hence $V = 0 = Z$. Thus $X = 0$, which implies that the orthogonal complement of
$[\cs,\cs]$ in $\cs_n \cap \cn$ is trivial. This establishes the proof of (\ref{derived-in-sn}).

Our next aim is to prove that
\begin{equation}\label{normal-in-v}
\nu_oF \subset \cv = \cg_\alpha\ .
\end{equation}
From (\ref{normal-in-n}) we already know that $\nu_oF \subset \cn$.
If $\nu_oF \cap \cv \neq 0$, we readily get $\nu_oF \subset \cv$, because
$\Ad(H_o)$ acts transitively on the unit sphere in $\nu_oF$ and preserves $\cv$.
Now assume that $\nu_oF \cap \cv = 0$. Then (\ref{derived-in-sn})
implies that the canonical projection of $[\cs,\cs] \subset \cv + \cz$ onto $\cv$ is
the entire space $\cv$. Thus, for each $V \in \cv$ there exists an element $V^\prime
\in \cz$ so that $V + V^\prime \in [\cs,\cs]$. Since $[\cs,\cs]$ is a subalgebra, we get
$$
[V,W] = [V + V^\prime, W + W^\prime] \in [\cs,\cs]
$$
for all $V,W \in \cv$. But since $[\cv,\cv] = [\cg_{\alpha},\cg_{\alpha}] = \cg_{2\alpha} = \cz$
this implies $\cz \subset [\cs,\cs]$ and hence $\nu_oF \subset \cv$. This establishes the proof
of (\ref{normal-in-v}).

From (\ref{normal-in-v}) we see that there exists a linear subspace $\cv_o$ of $\cv$ so that
$\cs_n = \ca + \cv_o + \cz$. Using the Lie algebra structure of $\ca + \cn$, we get:
\begin{equation}\label{sn-subalgebra}
\cs_n\ {\rm is\ a\ subalgebra\ of}\ \ca + \cn\ .
\end{equation}

Let $S_n$ be the connected subgroup of $AN$ with Lie algebra $\cs_n$. Our next aim is to show
that the orbit $S_n \cdot o$ of $S_n$ through $o$ coincides with the singular orbit $F$.
For this purpose we define
$$
\ct = \ck_x \cap \cs\ \subset \cs_c\ \ {\rm and}\ \ \ \cs^\prime = \R(\tilde{B} + B) + [\cs,\cs]\ ,
$$
where $B \in \ca$ and $\tilde{B} \in \cs_c$ are defined as above.
Since $\ct \subset \cs_c$, $\tilde{B} \in \cs_c$
and $\cs_c$ is abelian we have $[\ct,\tilde{B}] = 0$. And since $\ct \subset \ck_x$ and $\ck_x$
centralizes $\ca$ we have $[\ct,B] = 0$. Clearly, we also have $[\ct,[\cs,\cs]] \subset [\cs,\cs]
\subset \cs^\prime$ since $\ct \subset \cs$ and $[\cs,\cs] \subset \cs$. Altogether this
implies $[\ct,\cs^\prime] \subset \cs^\prime$. Moreover, since $\cs^\prime \subset \cs$, we
have $[\cs^\prime,\cs^\prime] \subset [\cs,\cs] \subset \cs^\prime$, which shows that
$\cs^\prime$ is a subalgebra of $\cs$. It follows that $\cs^\prime$ is an ideal in $\cs$ and
$\cs = \ct + \cs^\prime$ (semidirect sum). Let $S^\prime$ be the connected subgroup of $S$ with
Lie algebra $\cs^\prime$. Since $\ct \subset \ch_o$ we see that $S^\prime$ acts transitively on
$F$, i.e., $S^\prime \cdot o = F$.
For all $V+Z \in [\cs,\cs] \subset \cn = \cv + \cz$ we have
\begin{equation}\label{formula}
[\tilde{B},V+Z] + V + 2Z = [\tilde{B},V+Z] + [B,V+Z] = [\tilde{B} + B,V+Z] \in \cs^\prime
\end{equation}
since $\cs^\prime$ is a subalgebra. But $A$ and $V+Z$ are in $\cs_n$ according to
(\ref{a-in-sn}) and (\ref{derived-in-sn}), and since $\cs_n$ is a
subalgebra by (\ref{sn-subalgebra}),
we have $V + 2Z = [B,V+Z] \in \cs_n \cap \cn = [\cs,\cs] \subset \cs^\prime$
by (\ref{derived-in-sn}). By (\ref{formula}) this implies $[\tilde{B},V+Z] \in \cs^\prime$.
But $\tilde{B} \in \cs_c \subset \ck_x$ and thus $\ad(\tilde{B})$
leaves $\cv$ and $\cz$ invariant, which implies that $[\tilde{B},V+Z] \in \cs^\prime \cap \cn
\subset [\cs,\cs] \subset \cs_n$. We thus have proved that $\cs_n$ is normalized by $\tilde{B}$,
i.e., $[\tilde{B},\cs_n] \subset \cs_n$. Let ${\rm Exp}$ be the Lie exponential map of $\cg$.
We now get
$$
F = S^\prime \cdot o \subset {\rm Exp}(\R \tilde{B})S_n \cdot o = S_n {\rm Exp}(\R \tilde{B}) \cdot o
= S_n \cdot o
$$
since ${\rm Exp}(\R \tilde{B})$ normalizes $S_n$ and ${\rm Exp}(\R \tilde{B}) \cdot o \subset
H_o \cdot o = o$. Finally, by (\ref{dimensions})
the dimensions of $F$ and $S_n$ coincide, and since both $F$ and $S_n \cdot o$ are
complete, we must have $F = S_n \cdot o$. We thus have proved:

\begin{thm}\label{mainresult}
Let $H$ be a connected subgroup of $G = I^o(\F H^n)$ that acts on $\F H^n$
with cohomogeneity one and with a non-totally geodesic singular orbit $F$.
Then there exists a unique point $x \in M(\infty)$ that is fixed under
the induced action of $H$ on $M(\infty)$.
Let $o \in F$, $K$ the isotropy group of $G$ at $o$, and $\cg = \ck + \ca + \cn$
the Iwasawa decomposition of $\cg$ that is induced by $o$ and $x$.
Then there exists a subalgebra
$\cs$ of $\ca + \cn$ of the form $\cs = \ca + \cv_o + \cz$ with some linear
subspace $\cv_o$ of $\cv$, so that $F$ is the orbit of the connected subgroup $S$ of $AN$
with Lie algebra $\cs$.
\end{thm}

\section{The classification}

In this section we discuss the classification of cohomogeneity one actions on
noncompact symmetric spaces of rank one up to orbit equivalence. Recall that such an action
has either no singular orbit or exactly one singular orbit.

\smallskip
{\sc No singular orbit.}
In \cite{BT1} it was shown that there exist only two such actions without a singular orbit.
The first one is given by the action of the nilpotent group $N$ in an Iwasawa decomposition
$G = KAN$ of $G = I^o(\F H^n)$,
and the orbits form a horosphere foliation. The second one is given by the subgroup $S$
of $AN$ with Lie algebra $\cs = \ca + \cv_o + \cz$, where $\cv_o$ is a linear subspace
of $\cv$ with codimension one. The corresponding foliation has exactly one minimal leaf
and has been investigated in detail in \cite{Be}. In case of $\R H^n$ the minimal leaf
is a totally geodesic $\R H^{n-1} \subset \R H^n$.

\smallskip
{\sc Totally geodesic singular orbit.}
The cohomogeneity one actions on $\F H^n$ with a totally geodesic singular
orbit $F$ are given by:\\
\indent $M = \R H^n : F \in \{\pt,\R H^1,\ldots,\R H^{n-2}\}$;\\
\indent $M = \C H^n : F \in \{\pt,\C H^1,\ldots,\C H^{n-1},\R H^n\}$;\\
\indent $M = \H H^n : F \in \{\pt,\H H^1,\ldots,\H H^{n-1},\C H^n\}$;\\
\indent $M = \O H^2 : F \in \{\pt,\O H^1,\H H^2\}$.\\
Here, $\pt$ is a point in $\F H^n$, and the corresponding cohomogeneity one action is
just the action of the isotropy group of $I^o(\F H^n)$ at that point. More details about
this can be found in \cite{BB}.

\smallskip
{\sc Non-totally geodesic singular orbit.}
We now come to the classification of cohomogeneity one actions with a non-totally geodesic
singular orbit $F$. We will use the same notation as in the previous section.
Let $H$ be the connected component of the group of
isometries of $M$ that leave $F$ invariant. By Theorem \ref{mainresult}
there exists a unique point $x \in M(\infty)$ that is fixed under the induced action
of $H$ on $M(\infty)$. Let $o \in F$, $K$ the isotropy group of $G$ at $o$,
and $\cg = \ck + \ca + \cn$
the Iwasawa decomposition of $\cg$ that is induced by $o$ and $x$.
Using again Theorem \ref{mainresult}, there exists a subalgebra
$\cs$ of $\ca + \cn$ of the form $\cs = \ca + \cv_o + \cz$ with some linear
subspace $\cv_o$ of $\cv$,
so that $F$ is the orbit of the connected subgroup $S$ of $AN$
with Lie algebra $\cs$.
From the construction it is clear that the identity component of $H_o$
coincides with the identity component $N_K^o(\cs)$ of the normalizer $N_K(\cs)$
of $\cs$ in $K$. In order that $H$ acts with cohomogeneity
one it is therefore necessary and sufficient that the action of $N_K^o(\cs)$
on the normal space $\nu_oF$ is transitive on the unit sphere in $\nu_oF$.
Note that $N_K^o(\cs) \subset K_x$.
Since all Iwasawa decompositions of $\cg$ are conjugate to each other
under an inner automorphism
of $\cg$, it therefore remains to classify all subalgebras $\cs$ of $\ca + \cn$ of the form
$\cs = \ca + \cv_o + \cz$ with some linear subspace $\cv_o$ of $\cv$ such that
$N_K^o(\cs)$ acts transitively on the unit sphere in $\cv_o^\perp$, the orthogonal
complement of $\cv_o$ in $\cv$. This proves the first part of the following Theorem:

\begin{thm}\label{congruency}
Let $\cg = \ck + \ca + \cn$ be the Iwasawa decomposition induced by
$o \in M$ and $x \in M(\infty)$.
\begin{enumerate}
\item[{\rm (i)}]
Let $\cv_o$ be a linear subspace of $\cv$
so that $\dim \cv_o^\perp \geq 2$ and $N^o_{K_x}(\cv_o)$ acts
transitively on the unit sphere in $\cv_o^{\perp}$.
Then the connected subgroup of $G$ with Lie algebra
$N^o_{\cck_x}(\cv_o) + \ca + \cv_o + \cz$
acts on $M$ with cohomogeneity one so that the orbit
through $o$ is singular.
Furthermore, every cohomogeneity one action on $M$ with a non-totally
geodesic singular orbit can be obtained in this way
up to orbit equivalence.
\item[{\rm (ii)}]
Let $\cv_o$ and $\cv^\prime_o$ be linear subspaces of $\cv$ as in (i),
and assume that the corresponding cohomogeneity one actions
have non-totally geodesic singular orbits.
Then, these actions are orbit equivalent if and only if
there exists an isometry $k \in K_x$ so that ${\rm Ad}(k)\cv_o
= \cv_o^\prime$.
\end{enumerate}
\end{thm}

{\it Proof.} It remains to prove part (ii).
The ``if''-part of the statement is obvious.
Conversely, assume that the two cohomogeneity one actions are orbit equivalent.
Then the corresponding singular orbits, say $S$ and $S^\prime$,
are congruent under an isometry $k$ of $M$.
We may assume that $k$ fixes $o$.
By construction, the normalizers $N_G(S)$ and $N_G(S^\prime)$ fix $x$,
the point at infinity that determines our Iwasawa decomposition.
Then $k$ must fix $x$ as well, since $k N_G(S) k^{-1} = N_G(S^\prime)$ and
$x$ is the unique fixed point in $M(\infty)$ of $N_G(S)$ and of $N_G(S^\prime)$.
Therefore we conclude that ${\rm Ad}(k)\cv_o = \cv_o^\prime$.
\qed

\medskip
We now discuss the four different hyperbolic spaces individually.

\smallskip
\fbox{$M = \R H^n$} It follows from the classification of isoparametric hypersurfaces
in $\R H^n$ by Cartan \cite{Ca} that there exist no such actions. Since a singular orbit
of a cohomogeneity one action is necessarily minimal, one can also apply a result by
Di Scala and Olmos \cite{DO} stating that every minimal homogeneous submanifold of $\R H^n$ is
totally geodesic.
The classification also follows easily from Theorem \ref{congruency}:
Assume there is a cohomogeneity one action on $\R H^n$
with a non-totally geodesic singular orbit $F$.
Theorem \ref{congruency} implies that the action is orbit equivalent to
the $H$-action induced from $\ch = N^o_{\cck_x}(\cv_o) + \ca + \cv_o$
for some suitable subspace $\cv_o$ of $\cv$.
But for such an $H$-action the orbit $F = H \cdot o$ is totally geodesic,
which is a contradiction.

\smallskip
\fbox{$M = \C H^n$}
In this case the K\"{a}hler structure on $\C H^n$ induces a complex structure $J$ on
$\cv$ so that $\cv$ is isomorphic to $\C^{n-1}$ as a complex vector subspace.
Let $\cv_o$ be a linear subspace of $\cv$ so that $\dim_{\R} \cv_o^\perp \geq 2$.
Recall that the K\"{a}hler angle of a nonzero vector $v \in \cv_o^\perp \subset \C^{n-1}$
is defined as the angle between $Jv$ and $\cv_o^\perp$. In order that $N_K^o(\cs)$ acts
transitively on the unit sphere in $\cv_o^\perp$ it is necessary that the K\"{a}hler
angle of $\cv_o^\perp$ does not depend on the choice of the unit vector in $\cv_o^\perp$.
We thus assume that for all nonzero vectors $v \in
\cv_o^\perp$ the K\"{a}hler angle is equal to some $\varphi \in [0,\pi/2]$.
In the special case that $\varphi = 0$, $\cv_o^\perp$ is a complex subspace of $\cv$,
and if $\varphi = \pi/2$ then $\cv_o^\perp$ is a real subspace of $\cv$.
The subspaces of complex vector spaces with constant K\"{a}hler angle
have been classified in \cite{BB}. For $\varphi = 0$ we just have the complex subspaces
and for $\varphi = \pi/2$ the real subspaces, and in both cases the congruence classes
(under the action of $K_x = U(n-1)$ on $\cv = \C^{n-1}$)
are parametrized by the complex resp.\ real dimension. For $\varphi \in (0,\pi/2)$ there exists
exactly one congruence class of subspaces with constant K\"{a}hler angle $\varphi$
for each dimension $0 < 2k \leq n-1$.
For any such subspace the resulting
action on $\C H^n$ is of cohomogeneity one and $F$ is a non-totally geodesic
singular orbit unless $\varphi = 0$ (then $F$ is a totally geodesic complex
submanifold). Using Theorem \ref{congruency} we therefore conclude:

\begin{thm}
The moduli space of all cohomogeneity one actions on $\C H^n$, $n \geq 2$,
with a non-totally geodesic singular orbit (up to orbit equivalence) is isomorphic to
the disjoint union
$$
\{2,\ldots,n-1\} \cup \left((0,\pi/2) \times \{2k \mid k \in \Z\ ,\ 0 < 2k < n\}\right)\ .
$$
The integer in $\{2,\ldots,n-1\}$ indicates the codimension of the singular orbit if the
normal spaces are real, and the integer in $\{2k \mid k \in \Z\ ,\ 0 < 2k < n\}$ indicates
the codimension of the singular orbit if the normal spaces have constant K\"{a}hler angle
$\varphi \in (0,\pi/2)$.
\end{thm}

\begin{cor}
Any singular orbit of a cohomogeneity one action on $\C H^2$ is
totally geodesic.
\end{cor}

Note that by this result we now have a complete classification of
the homogeneous hypersurfaces in $\C H^n$ for all $n \geq 2$. In
view of \cite{BSF}, we call a submanifold normally homogeneous if
it is homogeneous and if the slice representation acts
transitively on the unit sphere in the normal bundle. A singular
orbit of a cohomogeneity one action is clearly a normally
homogeneous submanifold. The above shows that for each $k \in
\{2,\ldots,n-1\}$ there exists, up to holomorphic congruence,
exactly one normally homogeneous submanifold $F_k$ of $\C H^n$
with real normal bundle of rank $k$, and for each $k \in
\{1,\ldots,[(n-1)/2]\}$ and each $\varphi \in (0,\pi/2)$ there
exists exactly one, up to holomorphic congruence, normally
homogeneous submanifold $F_{k,\varphi}$ of $\C H^n$ with normal
bundle of rank $2k$ and constant K\"{a}hler angle $\varphi$.

\begin{thm}
Let $M$ be a homogeneous hypersurface in $\C H^n$, $n \geq 2$.
Then $M$ is holomorphically congruent to one of the following
hypersurfaces:
\begin{itemize}
\item[(1)] a tube of radius $r \in \R_+$ around the totally
geodesic $\C H^k \subset \C H^n$ for some $k \in
\{0,\ldots,n-1\}$; \item[(2)] a tube of radius $r \in \R_+$ around
the totally geodesic $\R H^n \subset \C H^n$; \item[(3)] a
horosphere in $\C H^n$; \item[(4)] the minimal ruled real
hypersurface $S$ determined by a horocycle in a totally geodesic
$\R H^2 \subset \C H^n$, or an equidistant hypersurface to $S$;
\item[(5)] a tube of radius $r \in \R_+$ around the normally
homogeneous submanifold $F_k$ of $\C H^n$ with real normal bundle
of rank $k$, $k \in \{2,\ldots,n-1\}$; \item[(6)] a tube of radius
$r \in \R_+$ around the normally homogeneous submanifold
$F_{k,\varphi}$ of $\C H^n$ with normal bundle of rank $2k \in
\{2,\ldots,2[(n-1)/2]\}$ and constant K\"{a}hler angle $\varphi
\in (0,\pi/2)$.
\end{itemize}
\end{thm}

\smallskip
\fbox{$M = \H H^n$}
In this case the quaternionic K\"{a}hler structure on $\H H^n$ induces a quaternionic structure $\cJ$ on
$\cv$ so that $\cv$ is isomorphic to $\H^{n-1}$ as a (right) quaternionic vector subspace.
Let $\cv_o$ be a linear subspace of $\cv$ so that $\dim_{\R} \cv_o^\perp \geq 2$.
In \cite{BB} the first author and Br\"uck
introduced the notion of a quaternionic K\"{a}hler angle, which
is defined as follows. Let $S^2$ be the two-sphere of all
almost Hermitian structures in $\cJ$.
For each nonzero vector $v \in \cv_o^\perp$ and each $J \in S^2$ denote by
$\varphi(v,J)$ the K\"{a}hler angle of $Jv$ and $\cv_o^\perp$ in the complex vector space
$(\cv,J)$. Since $S^2$ is compact, there exist minimum and maximum for these K\"{a}hler angles.
It was shown in \cite{BB} that for each nonzero $v$ there always exist a canonical basis $J_1,J_2,J_3$ of
elements in $S^2$ (i.e., $J_{\nu}J_{\nu+1} = J_{\nu+2} = - J_{\nu+1}J_{\nu}$, index modulo $3$)
such that $\varphi(v,J_1)$ is the minimum $\varphi_1(v)$ of these K\"{a}hler angles and $\varphi(v,J_3)$ is
the maximum $\varphi_3(v)$ of these K\"{a}hler angles. For any canonical basis with this
property the K\"{a}hler angle $\varphi_2(v) = \varphi(v,J_2)$ attains the same value. The
triple $\Phi(v) = (\varphi_1(v),\varphi_2(v),\varphi_3(v))$ of K\"{a}hler angles is called
the quaternionic K\"{a}hler angle of $\cv_o^\perp$ with respect to $v$. For a cohomogeneity
one action the quaternionic K\"{a}hler angle of $\cv_o^\perp$ must be independent of the
choice of the unit vector in $\cv_o^\perp$.
In \cite{BB} several examples of subspaces
of $\H^{n-1}$ with constant quaternionic K\"{a}hler angle were given, but a complete
classification is still missing. The examples are as follows:

\medskip
(a) $\Phi = (0,0,0)$.
The linear subspaces of $\cv$ with constant
quaternionic K\"ahler angle $\Phi =
(0,0,0)$ are the quaternionic subspaces. A linear subspace $V \subset \cv$ is quaternionic
if $JV \subset V$ holds for all $J \in \cJ$.
For each integer $k$ with $0 < k < n$ there exists exactly one (up to orbit
equivalence) cohomogeneity one action on $\H H^n$ with a
singular orbit $F$ of real codimension $4k$ with the property that the
normal spaces of $F$ have constant quaternionic K\"{a}hler angle
$\Phi = (0,0,0)$, and $F$ is congruent to the totally
geodesic $\H H^{n-k} \subset \H H^n$.

\medskip
(b) $\Phi = (0,\pi/2,\pi/2)$.
The linear subspaces of $\cv$ with constant
quaternionic K\"ahler angle $\Phi =
(0,\pi/2,\pi/2)$ are the totally complex subspaces.
A linear subspace $V \subset \cv$ is totally complex if there exists
an almost Hermitian structure $J_1 \in \cJ$ such that $J_1V \subset V$
and $JV \subset V^\perp$ for all $J \in \cJ$ perpendicular to $J_1$.
For each integer $k \in \{1,\ldots,n-1\}$ there exists exactly one (up to orbit
equivalence) cohomogeneity one action on $\H H^n$ with a non-totally
geodesic singular orbit $F$ of real codimension $2k$ with the property that the
normal spaces of $F$ have constant quaternionic K\"{a}hler angle
$\Phi = (0,\pi/2,\pi/2)$.

\medskip
(c) $\Phi = (\pi/2,\pi/2,\pi/2)$.
The linear subspaces of $\cv$ with constant quaternionic K\"ahler angle
$\Phi = (\pi/2,\pi/2,\pi/2)$ are the totally real subspaces.
A linear subspace $V \subset \cv$ is totally real
if $JV \subset V^\perp$ holds for all $J \in \cJ$.
For each integer $k \in \{2,\ldots,n-1\}$ there exists exactly one (up to orbit
equivalence) cohomogeneity one action on $\H H^n$ with a non-totally
geodesic singular orbit $F$ of real codimension $k$ with the property that the
normal spaces of $F$ have constant quaternionic K\"{a}hler angle
$\Phi = (\pi/2,\pi/2,\pi/2)$.

\medskip
(d) $\Phi = (0,0,\pi/2)$.
The linear subspaces of $\cv$ with constant quaternionic K\"ahler angle
$\Phi = (0,0,\pi/2)$ are the $3$-dimensional subspaces of the form
$(\Im \H)v$ for some unit vector $v \in \cv$.
There exists exactly one (up to orbit
equivalence) cohomogeneity one action on $\H H^n$ with a non-totally
geodesic singular orbit $F$ of real codimension $3$ with the property that the
normal spaces of $F$ have constant quaternionic K\"{a}hler angle
$\Phi = (0,0,\pi/2)$.

\medskip
(e) $\Phi = (\varphi,\pi/2,\pi/2)$, $\varphi \in (0,\pi/2)$. The
linear subspaces of $\cv$ with constant quaternionic K\"ahler
angle $\Phi = (\varphi,\pi/2,\pi/2)$, $\varphi \in (0,\pi/2)$, are
the linear subspaces with constant K\"{a}hler angle $\varphi$ in a
totally complex subspace $V$ of $\cv$. Here, the K\"{a}hler angle
in $V$ is measured with respect to the almost Hermitian structure
$J_1$ as described in (b). For each integer $k \in
\{1,\ldots,[(n-1)/2]\}$ and each $\varphi \in (0,\pi/2)$ there
exists exactly one (up to orbit equivalence) cohomogeneity one
action on $\H H^n$ with a non-totally geodesic singular orbit $F$
of real codimension $2k$ with the property that the normal spaces
of $F$ have constant quaternionic K\"{a}hler angle $\Phi =
(\varphi,\pi/2,\pi/2)$.

\medskip
(f) $\Phi = (0,\varphi,\varphi)$, $\varphi \in (0,\pi/2)$. The
linear subspaces of $\cv$ with constant quaternionic K\"ahler
angle $\Phi = (0,\varphi,\varphi)$, $\varphi \in (0,\pi/2)$, are
the complexifications of linear subspaces with constant K\"{a}hler
angle $\varphi$ in a totally complex subspace $\cw$ of $\cv$. More
precisely, let $J_2 \in \cJ$ be an almost Hermitian structure and
consider $\cv$ as the complexification of $\cw$ with respect to an
almost Hermitian structure $J_1 \in \cJ$ orthogonal to $J_2$, that
is $\cv = \cw + J_1 \cw$ with a $J_2$-invariant linear subspace
$\cw \subset \cv$. Let $W$ be a linear subspace of the complex
vector space $(\cw,J_2)$ with constant K\"{a}hler angle $\varphi$.
Then the complexification of $W$ with respect to $J_1$ is a linear
subspace of $\cv$ with constant quaternionic K\"ahler angle $\Phi
= (0,\varphi,\varphi)$. For each integer $k \in
\{1,\ldots,[(n-1)/2]\}$ and each $\varphi \in (0,\pi/2)$ there
exists exactly one (up to orbit equivalence) cohomogeneity one
action on $\H H^n$ with a non-totally geodesic singular orbit $F$
of real codimension $4k$ with the property that the normal spaces
of $F$ have constant quaternionic K\"{a}hler angle $\Phi =
(0,\varphi,\varphi)$.

\medskip
We conjecture that each cohomogeneity one action on $\H H^n$ with
a non-totally geodesic singular orbit is orbit equivalent to one of
these examples.
This is true for $n=2$, and for the case that the singular orbit has
codimension $2$.

\begin{thm}
The moduli space of all cohomogeneity one actions on $\H H^2$
with a non-totally geodesic singular orbit (up to orbit equivalence) is isomorphic to
the set $\{2,3\}$.
The number $k \in \{2,3\}$ parametrizes the unique (up to
orbit equivalence) cohomogeneity one
action on $\H H^2$ with a non-totally geodesic
singular orbit of codimension $k$.
\end{thm}

{\it Proof.} In case of $\H H^2$ the quaternionic vector space $\cv$ has
quaternionic dimension one. It is easy to see that every $2$-dimensional
subspace of a one-dimensional quaternionic subspace has constant
quaternionic K\"{a}hler angle $\Phi = (0,\pi/2,\pi/2)$, and every
$3$-dimensional
subspace of a one-dimensional quaternionic subspace has constant
quaternionic K\"{a}hler angle $\Phi = (0,0,\pi/2)$. The result then
follows from (b) and (d) above. Note that codimension $4$ occurs for
quaternionic K\"{a}hler angle $\Phi = (0,0,0)$, which leads to a totally
geodesic singular orbit.
\qed

\begin{thm}
The moduli space of all cohomogeneity one actions on $\H H^n$, $n > 2$,
with a non-totally geodesic singular orbit with codimension $2$
(up to orbit equivalence) is isomorphic to the closed interval $[0,\pi/2]$.
The number $\varphi \in [0,\pi/2]$ parametrizes the unique (up to
orbit equivalence) cohomogeneity one
action on $\H H^n$ with a non-totally geodesic
singular orbit of codimension $2$ for which
the normal spaces have constant quaternionic K\"{a}hler angle
$\Phi = (\varphi,\pi/2,\pi/2)$.
\end{thm}

{\it Proof.} Every $2$-dimensional subspace of $\cv$ has
constant quaternionic K\"{a}hler angle
$\Phi = (\varphi,\pi/2,\pi/2)$ for some $\varphi \in [0,\pi/2]$. The result
then follows from (b), (c) and (e) above. \qed

\smallskip
\fbox{$M = \O H^2$}
In \cite{BB} the first author and Br\"uck
classified all subspaces $\cv_o$ of $\cv = \R^8$ for which there exists a
subgroup of $K_x = \Spin(7)$ that acts transitively on the unit
sphere in $\cv_o^\perp$. In fact, any subspace $\cv_o$ of $\cv$
with dimension $k \in \{1,2,4,5,6\}$ has this property, but there
are no $3$-dimensional subspaces with this property.
We denote by $G_k^+(\R^8)$ the Grassmann manifold of oriented $k$-planes
in $\R^8$, and by $G_k(\R^8)$ the Grassmann manifold of (unoriented)
$k$-planes in $\R^8$. It is clear that $G_k^+(\R^8)$ is a two-fold
covering of $G_k(\R^8)$, and that there is a natural isomorphism between
the Grassmann manifolds of $k$- and $(8-k)$-planes.

The Lie group $\Spin(7)$ acts on $\R^8$ by its irreducible $8$-dimensional
spin representation. This naturally induces actions of $\Spin(7)$ on
$G_k^+(\R^8)$ and $G_k(\R^8)$. For $k=1$,
it was proved by Borel \cite{Bo} that
$\Spin(7)$ acts transitively on $G_1^+(\R^8) = S^7$ and that $S^7 =
\Spin(7)/G_2$. For $k=2$ we also have a transitive action,
so that $G_2^+(\R^8) = \Spin(7)/U(3)$ (see e.g.\
\cite{Br}), and hence also $G_6^+(\R^8) = \Spin(7)/U(3)$.
Also for $k=3$ the action is transitive, and we have
$G_3^+(\R^8) = \Spin(7)/SO(4)$ (see e.g.\ \cite{On}), and thus also
$G_5^+(\R^8) = \Spin(7)/SO(4)$.

The action of $\Spin(7)$ on $G_4^+(\R^8)$ is not transitive,
but of cohomogeneity one (see \cite{Br} and \cite{HaL} for details).
One singular orbit of this action consists of the so-called
Cayley 4-planes in $\O$ introduced by Harvey and Lawson \cite{HaL}.
The submanifolds of $\O$ all of whose tangent spaces are Cayley 4-planes
are so-called Cayley submanifolds of $\O$ and provide a beautiful
example of a calibrated geometry. This singular orbit is isomorphic
to $\Spin(7)/(SU(2)^3/\Z_2)$, and the second singular orbit consists just
of the Cayley 4-planes with opposite orientation. This can also be seen in
the following way. Let $V \in G_3^+(\R^8)$ be an oriented $3$-plane in $\R^8$.
We know from the above that $Spin(7)$ acts transitively on $G_3^+(\R^8)$ and
the isotropy group at $V$ is some $SO(4) \subset Spin(7)$. There is a
unique unit vector $\xi$ in the orthogonal complement $V^\perp$ of $V$ in
$\R^8$ so that the $4$-plane $V \oplus \R\xi$ is a Cayley 4-plane. Then
$V \oplus \R(-\xi)$ is the same $4$-plane with opposite orientation.
The action of $SO(4)$ on the unit sphere $S^4$ in $V^\perp$ is the
standard action determined by the two fixed points $\pm\xi$. The principal
orbits are the $3$-spheres in $S^4$ with center $\xi$. Each such orbit parametrizes
in a canonical way a set of oriented $4$-planes in $\R^8$ containing the
$3$-dimensional subspace $V$. We now turn to the induced action of $Spin(7)$ on the
Grassmannian $G_4(\R^8)$ of unoriented $4$-planes in $\R^8$.
This action is clearly of cohomogeneity one as well. The two
singular orbits on $G_4^+(\R^8)$ become identified under the two-fold
covering map $G_4^+(\R^8) \to G_4(\R^8)$, and provide one singular orbit
of the action. The second singular orbit in $G_4(\R^8)$
is the projection of the principal orbit on $G_4^+(\R^8)$
containing $4$-planes of both orientations. This orbit contains the $4$-planes that
are constructed from the unique totally geodesic principal orbit of the
$SO(4)$-action on $S^4 \subset V^\perp$. The second singular orbit in
$G_4(\R^8)$ is therefore a $2$-fold subcovering of any principal orbit
and thus has the same dimension as the principal orbits.

From Theorem \ref{congruency}
it is clear that if the action of $\Spin(7)$ is transitive on
$G_k(\R^8)$, then all cohomogeneity one actions constructed from a
$k$-dimensional subspace of $\cv$ are orbit equivalent. In the case
$k=4$, the cohomogeneity one actions induced from a $4$-dimensional
subspace of $\cv$ up to orbit equivalence are in one-to-one
correspondence with the orbits of the action of $\Spin(7)$ on $G_4(\R^8)$.
Altogether this now implies:

\begin{thm}
The moduli space of all cohomogeneity one actions on $\O H^2$
with a non-totally geodesic singular orbit (up to orbit equivalence) is isomorphic to
the disjoint union
$$
\{2,3,6,7\} \cup (\{4\}\times [0,1])\ .
$$
The number $k \in \{2,3,6,7\}$ parametrizes the unique (up to
orbit equivalence) cohomogeneity one
action on $\O H^2$ with a singular orbit of codimension $k$.
The set $\{4\} \times [0,1]$ parametrizes
the cohomogeneity one actions on $\O H^2$ with a
singular orbit of codimension $4$ (up to orbit equivalence).
\end{thm}

The above result says that for each $k \in \{2,3,6,7\}$ there
exists exactly one, up to isometric congruence, normally
homogeneous submanifold $F_k$ of $\O H^2$ with normal bundle of
rank $k$, and for each $\varphi \in [0,1]$ there exists exactly
one, up to isometric congruence, normally homogeneous submanifold
$F_{4,\varphi}$ of $\O H^2$ with normal bundle of rank $4$. We now
have a complete classification of the homogeneous hypersurfaces in
the Cayley hyperbolic plane.

\begin{thm}
Let $M$ be a homogeneous hypersurface in $\O H^2$. Then $M$ is
isometrically congruent to one of the following hypersurfaces:
\begin{itemize}
\item[(1)] a geodesic hypersphere of radius $r \in \R_+$ in $\O
H^2$; \item[(2)] a tube of radius $r \in \R_+$ around the totally
geodesic $\O H^1 \subset \O H^2$; \item[(3)] a tube of radius $r
\in \R_+$ around the totally geodesic $\H H^2 \subset \O H^2$;
\item[(4)] a horosphere in $\O H^2$; \item[(5)] the minimal
homogeneous hypersurface $S$ in $\O H^2$, or an equidistant
hypersurface to $S$; \item[(6)] a tube of radius $r \in \R_+$
around the normally homogeneous submanifold $F_k$ of $\O H^2$ with
normal bundle of rank $k$, $k \in \{2,3,6,7\}$; \item[(7)] a tube
of radius $r \in \R_+$ around the normally homogeneous submanifold
$F_{4,\varphi}$ of $\O H^2$ with normal bundle of rank $4$ and
$\varphi \in [0,1]$.
\end{itemize}
\end{thm}

\bigskip
\noindent {\small Department of Mathematics, University College, Cork, Ireland, email: j.berndt@ucc.ie}

\bigskip
\noindent {\small Department of Mathematics, Hiroshima University,
1-3-1 Kagamiyama, Higashi-Hiroshima, 739-8526, Japan,
email: tamaru@math.sci.hiroshima-u.ac.jp}

\end{document}